\documentclass[a4paper]{article}
\usepackage{a4, ams, amsmath, amssymb, graphics, subfigure, fullpage, color}%
\usepackage[pdftex]{graphicx}
\usepackage[colorlinks]{hyperref}
\hypersetup{colorlinks,citecolor=green,filecolor=black,linkcolor=blue,urlcolor=blue}
\usepackage{fancyhdr}
\usepackage{multirow}
\usepackage{amsthm}

\newtheorem{lemma}{Lemma}

\usepackage{lineno}
\usepackage{upgreek} 
\usepackage{lscape}

\newcommand{\cN}{\mathcal{N}}

\DeclareMathOperator{\e}{\textrm{e}}

\newcommand{\beq}{\begin{equation}}
\newcommand{\eeq}{\end{equation}}
\newcommand{\beqn}{\begin{eqnarray}}
\newcommand{\eeqn}{\end{eqnarray}}
\newcommand{\bfig}{\begin{figure}}
\newcommand{\efig}{\end{figure}}
\newcommand{\btab}{\begin{table}}
\newcommand{\etab}{\end{table}}

\ifdefined \up 
	\renewcommand{\up}{\textrm{up}}
\else 
	\DeclareMathOperator{\up}{\textrm{up}} 
\fi

\usepackage{lscape}
\usepackage{lineno}
\usepackage[utf8]{inputenc}
\usepackage[T1]{fontenc}

\usepackage{natbib}

\newcommand{\MyBib}{./ConditionalSamplingR2}

\title{{\bf {\Large Closed form expression of the multivariate standard Normal distribution under a weighted sum constraint}}}
\author{Fr\'ed\'eric Vrins\thanks{Contact information: Voie du Roman Pays 34, B-1348 Louvain-la-Neuve, Belgium. E-mail: \href{mailto:frederic.vrins@uclouvain.be}{frederic.vrins@uclouvain.be} The author thanks Damiano Brigo for interesting discussions around this question.}\vspace{0.2cm}\\ Louvain Finance Center \& CORE\\ Universit\'e catholique de Louvain.}
\date{\today}

\begin{document}
\maketitle
\begin{abstract} In this letter we derive the $(n-1)$-dimensional distribution corresponding to a $n$-dimensional i.i.d. Normal standard vector $Z=(Z_1,Z_2,\ldots,Z_n)$ subjected to the weighted sum constraint $\sum_{i=1}^n w_i Z_i=c$, $w_i\neq 0$. We first address the $n=2$ case before proceeding with the general $n\geq 2$ case. The resulting distribution is a Normal distribution whose mean vector $\mu$ and covariance matrix $\Sigma$ are explicitly derived as a function of $w_1,\ldots,w_n,c$. The derivation of the density relies on a very specific positive definite matrix for which the determinant and inverse can be computed analytically.
\end{abstract}

\section{Introduction}

Factor models are extensively used in statistical modeling. In banking and finance for instance, it is a standard procedure to introduce a dependence structure among loans in credit risk modeling, see e.g. Li's model~\cite{Li00} but also \cite{Hull04,Ande04,Vrins09c,Laur16}, just to name a few. 
In such models, the credit worthiness of the $i$-th entity is typically modeled as a random variable $X_i$ defined as a weighted sum of common factors $(Y_1,\ldots,Y_J)$ accounting for the state of the global economy, the sector, the region, etc, and an idiosynchratic variable $\epsilon_i$. In the popular case of a Gaussian copula model, all these factors are Normally distributed. The $Y_j$ factors do not need to be independent, but can be decomposed (via a Cholesky transform) as a weighted sum of independent Normal risk factors $\tilde{Z}:=(\tilde{Z}_1,\tilde{Z}_2,\ldots,\tilde{Z}_n)$. 
%
The knowledge of a default event of the $i$-th reference entity reveals that the credit worthiness variable $X_i$ reached the (assumed to be known) default threshold $c_i$. The event $X_i=c_i$ carries some information about the distribution of the underlying factors \textit{in that specific state}. In particular, the vector $\tilde{Z}$ is no longer standard Normal \textit{being told that $X_i=c_i$}; risk measures (e.g. value-at-risk) of the portfolio built from the outstanding loans might be strongly impacted by this information. This raises the following question: given a value $y$ for the weighted sum $w'Z$, what is the distribution of $Z$ ? 
Even if the analytical form of the conditional distribution is unknown, it is of course straightforward to sample a vector $Z$ of $n$ Normal variables such that the weighted sum is $y$. One possibility is to sample $Z_j\sim\cN(0,1)$ for $j\in\{1,2,\ldots,n-1\}$ and then set $Z_n=(y-\sum_{i=1}^{n-1}w_i Z_i)/w_n$. Another possibility would be to sample a vector of $n$ i.i.d. standard Normal variables $\tilde{Z}=(\tilde{Z_1},\tilde{Z_2},\ldots,\tilde{Z_n})$, compute $\tilde{y}:=w'\tilde{Z}$ and rescale the $\tilde{Z}$ to set $Z=\frac{y}{\tilde{y}}\tilde{Z}$. Alternatively, one could take $Z_i=\tilde{Z}_i+(y-\tilde{y})/(nw_i)$. However, none of these approaches yield the correct answer. The later requires the knowledge of the conditional distribution.\medskip 

In this letter, we derive the conditional distribution associated to the $(w'Z=c)$-\textit{slice} of the $n$-dimensional standard Normal density when $w_i\neq 0$ for all $i\in\{1,2,\ldots,n\}$. Interestingly, it is a $(n-1)$-Normal whose mean vector and covariance matrix can be computed in closed form, respectively given by



$$\mu(c,w):=\frac{c}{w_n^2\|w\|_2^2}\hbox{diag}(ww')~~,~~\Sigma_{i,j}(w):=\frac{w_i^2}{\|w\|_2^2} \left(\delta_{ij}\left(\|w\|_2^2-w_i^2\right)+(\delta_{ij}-1)w_j^2\right)\; .$$


The distribution of $Z$ can be obtained by simple rescaling of that of $X$ as $X=DZ$ where $D=\hbox{diag}(w)$ is an invertible diagonal matrix. We first address the $(n=2)$-case before moving to the general case $n\geq 2$. The result derives from the analytical properties of a square positive definite matrix having a very specific form.

\section{Bivariate case}

We are looking for the distribution of $(X_1,X_2)$ given that $X_1+X_2=c$. If $Z_i\sim\cN(0,1)$ iid, then $X_i\sim\cN(1,w_i)$ are independent normal variables with standard deviation $w_i$. We note $\phi(x;\mu,\sigma)$ the density associated to $\cN(\mu,\sigma)$. \medskip

We first compute the conditional density using Bayes

$$f_{X}(x|x_1+x_2=c)=f_{X_1,X_2}(x_1,x_2|x_1+x_2=c)=\frac{f_{X_1,X_2}(x_1,x_2;x_1+x_2=c)}{f_{x_1+x_2}(c)}$$

where the denominator is the centered Normal density with standard deviation $\sqrt{w_1^2+w_2^2}$:

$$k_1(c,w):=f_{x_1+x_2}(c)=\phi\left(c;0,\sqrt{w_1^2+w_2^2}\right)\;.$$

The numerator reads

$$\frac{1}{\sqrt{2\pi} w_1}\e^{-\frac{x_1^2}{2 w_1^2}}\frac{1}{\sqrt{2\pi} w_2}\e^{-\frac{x_2^2}{2 w_2^2}}=\frac{1}{2\pi w_1w_2}\e^{-\frac{(x_1/w_1)^2+((c-x_1)/w_2)^2}{2}}\;.$$

One can thus develop and complete the square to get

$$\frac{1}{2\pi w_1w_2}\e^{-\frac{x_1^2}{2w_1^2}}\e^{-\frac{c^2-2x_1c+x_1^2}{2w_2^2}}=\frac{\e^{-\frac{c^2}{2w_2^2}}}{2\pi w_1w_2}\e^{-\frac{x_1^2}{2}\left(\frac{1}{w_1^2}+\frac{1}{w_2^2}\right)+\frac{x_1c}{w_2^2}}$$

and

$$\frac{\e^{-\frac{c^2}{2w_2^2}}\e^{\frac{c^2}{2w_2^4\left(\frac{1}{w_1^2}+\frac{1}{w_2^2}\right)}}}{2 w_1w_2}\e^{-\frac{1}{2}\left(\frac{1}{w_1^2}+\frac{1}{w_2^2}\right)\left(x_1-\frac{c}{w_2^2\left(\frac{1}{w_1^2}+\frac{1}{w_2^2}\right)}\right)^2}=k_1(c,w)\phi\left(x,\frac{c}{w_2^2\left(\frac{1}{w_1^2}+\frac{1}{w_2^2}\right)},\left(\sqrt{\frac{1}{w_1^2}+\frac{1}{w_2^2}}\right)^{-1}\right)\;.$$

%


Hence, the conditional density $f(x_1,x_2|c)$ of $(X_1,X_2)$ at $(x_1,c-x_1)$ is given by $f(x_1)$ where


\beqn
f(x)&:=&\phi\left(x;\mu(x,w),\sigma(w)\right)\nonumber\\
\sigma(w)&:=&\left(\sqrt{\frac{1}{w_1^2}+\frac{1}{w_2^2}}\right)^{-1}\nonumber\\
\mu(c,w)&:=&\frac{c}{w_2^2\left(\frac{1}{w_1^2}+\frac{1}{w_2^2}\right)}=\frac{c}{w_2^2}\sigma^2(w)\nonumber\;.
\eeqn

%

\section{Extension to higher dimensions}

As before we compute the conditional density starting from Bayes' theorem,
$$f_{X}\left(x\left|\sum_{i=1}^n x_i=c\right.\right):=f_{X_1,\ldots,X_n}\left(x_1,\ldots,x_n\left|\sum_{i=1}^n x_i=c\right.\right)=\frac{f_{X_1,\ldots,X_n}\left(x_1,\ldots,x_n;\sum_{i=1}^n x_i=c\right)}{f_{\sum_{i=1}^n x_i}(c)}\;.$$

The denominator collapses to the one-dimensional centered Normal density with variance $w'w$:

$$k_1(c,w):=f_{\sum_{i=1}^n x_i}(c)=\phi\left(c;0,\sqrt{\sum_{i=1}^n w_i^2}\right)\;.$$

The numerator can be written as

$$\left(\prod_{i=1}^{n-1}\frac{\e^{-\frac{x_i^2}{2 w_i^2}}}{\sqrt{2\pi} w_i}\right)\frac{\e^{-\frac{\left(c-\sum_{i=1}^{n-1}x_i\right)^2}{2 w_n^2}}}{\sqrt{2\pi} w_n}=\left(\prod_{i=1}^{n-1}\frac{\e^{-\frac{x_i^2}{2 w_i^2}}}{\sqrt{2\pi} w_i}\right)\frac{\e^{-\frac{c^2-2c\sum_{i=1}^{n-1}x_i+\sum_{i=1}^{n-1}x_i^2+\sum_{i=1}^{n-1}\sum_{j=1,j\neq i}^{n-1}x_ix_j}{2 w_n^2}}}{\sqrt{2\pi} w_n}$$
$$=\left(\prod_{i=1}^{n}\frac{1}{\sqrt{2\pi} w_i}\right)\exp\left\{-\frac{1}{2}\sum_{i=1}^{n-1}\left(\left(\frac{1}{w_i^2}+\frac{1}{w_n^2}\right)x_i^2+\frac{x_i}{w_n^2}\sum_{j=1,j\neq i}^{n-1}x_j-\frac{2c}{w_n^2}x_i\right)-\frac{c^2}{2w_n^2}\right\}\;.$$

Hence, the conditional density looks like that of a $(n-1)$-th dimensional Normal pdf:

\beq
f_{X}\left(x\left|\sum_{i=1}^n x_i=c\right.\right)=k_2\exp\left\{-\frac{1}{2}\sum_{i=1}^{n-1}\left(\left(\frac{1}{w_i^2}+\frac{1}{w_n^2}\right)x_i^2+\frac{x_i}{w_n^2}\sum_{j=1,j\neq i}^{n-1}x_j-\frac{2c}{w_n^2}x_i\right)-\frac{c^2}{2w_n^2}\right\}\label{eq:condpdf1}
\eeq

where

$$k_2:=k(c,w)\e^{\frac{c^2}{2\sum_{i=1}^n w_i^2}}~~\hbox{ and }~~k(c,w):=\frac{1}{(\sqrt{2\pi})^{n-1}}\frac{\sqrt{\sum_{i=1}^n w_i^2}}{\prod_{i=1}^{n} w_i}\; .$$

In order for this density to belong to the Normal family, it needs to take the form of $\phi(x;\mu,\Sigma)$ where $\Sigma$ is a valid (positive definite) covariance matrix. In the sequel, we prove that $f_{X}\left(x\left|\sum_{i=1}^n x_i=c\right.\right)$ does indeed have such a form and confirm that the corresponding matrix $\Sigma$ is positive definite by determining the entries $\alpha_{i,j}$ of $\Sigma^{-1}$, the inverse of the $(n-1)$-dimensional covariance matrix $\Sigma$, and showing that $\Sigma^{-1}$ is invertible and positive definite. Moreover, we compute analytically $\Sigma$ and its determinant $|\Sigma|$ as well as the corresponding mean vector $\mu=(\mu_1,\mu_2,\ldots,\mu_{n,-1})$.\medskip

We start with the development of the Normal density of dimension $n-1$ :

\beqn
\phi(x;\mu,\Sigma)
&=&K\exp\left\{-\frac{1}{2}\left(\sum_{i=1}^{n-1}\sum_{j=1}^{n-1}\alpha_{i,j}x_ix_j-\sum_{i=1}^{n-1}\sum_{j=1}^{n-1}\alpha_{i,j}\mu_ix_j-\sum_{i=1}^{n-1}\sum_{j=1}^{n-1}\alpha_{i,j}\mu_jx_i+\sum_{i=1}^{n-1}\sum_{j=1}^{n-1}\alpha_{i,j}\mu_i\mu_j\right)\right\}\nonumber\\
&=&K\exp\left\{-\frac{1}{2}\left(\sum_{i=1}^{n-1}\sum_{j=1}^{n-1}\alpha_{i,j}x_ix_j-\sum_{i=1}^{n-1}\sum_{j=1}^{n-1}\alpha_{j,i}\mu_jx_i-\sum_{i=1}^{n-1}\sum_{j=1}^{n-1}\alpha_{i,j}\mu_jx_i+\sum_{i=1}^{n-1}\sum_{j=1}^{n-1}\alpha_{i,j}\mu_i\mu_j\right)\right\}\nonumber\\
&=&K\exp\left\{-\frac{1}{2}\sum_{i=1}^{n-1}\left(\alpha_{i,i}x_i^2+x_i\sum_{j=1,j\neq i}^{n-1}\alpha_{i,j}x_j-x_i\sum_{j=1}^{n-1}(\alpha_{j,i}+\alpha_{i,j})\mu_j+\mu_i\sum_{j=1}^{n-1}\alpha_{i,j}\mu_j\right)\right\}\label{eq:jpdf1}
\eeqn

where $K:=1/\sqrt{(2\pi)^{n-1} |\Sigma|}$. To determine the expression of the covariance matrix and mean vector of the conditional density \eqref{eq:condpdf1} (assuming it is indeed Normal), it remains to determine the entries of $\mu,\Sigma^{-1}$ by inspection, comparing the expression of conditional density in \eqref{eq:condpdf1} with that of the multivariate Normal~\eqref{eq:jpdf1}.\medskip

Leaving only $k(c,w)$ as a factor in front of the exponential in \eqref{eq:condpdf1}, the independent term (i.e. the term that does not appear as a factor of any $x_i$)  reads without any loss of generality as

$$\frac{c^2}{2\sum_{i=1}^n w_i^2}-\frac{c^2}{2 w_n^2}=-\frac{c^2}{2w_n^2}\frac{\sum_{i=1}^{n-1} w_i^2}{\sum_{i=1}^{n} w_i^2}=-\frac{c^2}{2w_n^2\sum_{j=1}^{n} w_j^2}\sum_{i=1}^{n-1} \gamma_i w_i^2$$

for any $(\gamma_1,\gamma_2,\ldots,\gamma_{n-1})$ satisfying $\sum_{i=1}^{n-1} \gamma_iw_i^2=\sum_{i=1}^{n-1} w_i^2$.\footnote{The constant case $\gamma_i=1$ might be a solution but it is not guaranteed at this stage.}
Comparing \eqref{eq:condpdf1} and~\eqref{eq:jpdf1}, it comes that the expression

\beq
\left(\frac{1}{w_i^2}+\frac{1}{w_n^2}\right)x_i^2+\frac{x_i}{w_n^2}\sum_{j=1,j\neq i}^{n-1}x_j-\frac{2c}{w_n^2}x_i+\frac{c^2\gamma_i w_i^2}{w_n^2\sum_{j=1}^{n} w_j^2}\label{eq:Sol1}
\eeq

must agree with

\beq
\alpha_{i,i}x_i^2+x_i\sum_{j=1,j\neq i}^{n-1}\alpha_{i,j}x_j-x_i\sum_{j=1}^{n-1}(\alpha_{i,j}+\alpha_{j,i})\mu_j+\mu_i\sum_{j=1}^{n-1}\alpha_{i,j}\mu_j\label{eq:Sol2}
\eeq

for all $x_1,x_2,\ldots,x_{n-1}$. Equating the $x_ix_j$ terms in~\eqref{eq:Sol1} and~\eqref{eq:Sol2} uniquely determines the components of $\Sigma^{-1}$, $\alpha_{i,i}:=(\Sigma^{-1})_{i,i}=\frac{1}{w_i^2}+\frac{1}{w_n^2}$ and $\alpha_{i,j\neq i}:=(\Sigma^{-1})_{i,j\neq i}=\frac{1}{w_n^2}$. It remains to show that $k(c,w)=K$, to find the expressions of the $\mu_i$'s from the $x_i$ terms, provide the expression of $\Sigma$ by inverting $\Sigma^{-1}$ and finally, to check that the independent terms in~\eqref{eq:Sol1} and~\eqref{eq:Sol2} agree and that the implied $\gamma_i$'s comply with $\sum_{i=1}^{n-1} \gamma_iw_i^2=\sum_{i=1}^{n-1} w_i^2$. To that end, we rely on the following lemma (proven in the end of the paper).

\begin{lemma}Let $\delta_{ij}$ be the Kronecker delta and $A(m)$ denote a matrix with $(i,j)$ elements $A_{ij}(m)=a_i\delta_{ij}+a_0$, $a_{k}>0$ for all $k\in\{0,1,\ldots,m\}$. Define $\pi(m):=\prod_{k=0}^{m} a_k$ and $s(m):=\sum_{k=0}^m 1/a_k$. Then:

\begin{itemize}
\item[(i)] $A(m)$ is positive definite ;
\item[(ii)] its determinant is given by

$$|A(m)|=\sum_{k=0}^{m} \prod_{j=0,j\neq k}^m a_j=\pi(m)s(m)\; ;$$

\item[(iii)] the elements of the inverse $B(m):=(A(m))^{-1}$ are given by

$$B_{i,j}(m)=\frac{1}{a_is(m)}\left(\delta_{ij}\frac{a_is(m)-1}{a_i}+\frac{\delta_{ij}-1}{a_j}\right)\;.$$


\end{itemize}
\end{lemma}

As $\Sigma^{-1}$ takes the form $A(n-1)$ with $a_0\leftarrow 1/w_n^2$ and $a_i\leftarrow 1/w_i^2$ for $i\in\{1,2,\ldots,n-1\}$ we can call Lemma 1 $(i)$ to show that $\Sigma^{-1}$ is symmetric and positive definite, proving that $\Sigma$ is a valid covariance matrix satisfying $|\Sigma|>0$. From Lemma 1 $(ii)$, $k(c,w)=K$ as\footnote{Observe that in $A(n-1)$ the summation and product indices agree with that of the $a_i$, i.e. range from 0 to $n-1$, but the index of $w_i$ ranges from $1$ to $n$.}

$$|\Sigma^{-1}|=\left(\prod_{j=1}^{n} \frac{1}{w_j^2}\right)\left(\sum_{i=1}^{n} w_i^2\right)=\frac{\sum_{i=1}^{n} w_i^2}{\prod_{j=1}^{n} w_j^2}~~\Rightarrow ~~1/\sqrt{|\Sigma|}=\sqrt{|\Sigma^{-1}|}= \frac{\sqrt{\sum_{k=1}^{n} w_k^2}}{\prod_{k=1}^{n} w_k}\; .$$

We can then use Lemma 1 $(iii)$ to determine $B(n-1)$, the elements $\beta_{i,j}$ of $\Sigma$. Setting $\|w\|_2:=\sqrt{\sum_{k=1}^{n} w_k^2}$,

\beqn
\beta_{i,j}
&=&\frac{w_i^2}{\|w\|_2^2} \left(\delta_{ij}\left(\|w\|_2^2-w_i^2\right)+(\delta_{ij}-1)w_j^2\right)\; .\nonumber
\eeqn

Finally, the mean vector is obtained by equating the $x_i$ terms in~\eqref{eq:Sol1} and~\eqref{eq:Sol2}. Using that $\Sigma^{-1}$ is symmetric, we observe that for all $i \in\{1,2,\ldots,n-1\}$:

\beq
\frac{2c}{w_n^2}=2\sum_{j=1}^{n-1}\alpha_{i,j}\mu_j\Rightarrow \sum_{j=1}^{n-1}\alpha_{i,j}\mu_j=\frac{c}{w_n^2}\;.\label{eq:sumamu}
\eeq
Hence, $\Sigma^{-1}\mu=\frac{c}{w_n^2}{\bf{1}}_{n-1}$ where ${\bf{1}}_m$ is the $m$-dimensional column vector with $m$ entries all set to 1 so that $\mu_i=\frac{c}{w_n^2}\sum_{j=1}^{n-1}\beta_{i,j}=\frac{cw_i^2}{\|w\|_2^2}$.\medskip

It remains to check that these expressions for $\mu$ and $\Sigma$ also comply with the independent term. Equating the independent terms of \eqref{eq:Sol1} and~\eqref{eq:Sol2} and calling ~\eqref{eq:sumamu} yields

$$\frac{c^2\gamma_i w_i^2}{w_n^2\|w\|_2^2}=\mu_i\frac{c}{w_n^2}~~\Rightarrow~~\mu_i=\frac{c\gamma_iw_i^2}{\|w\|_2^2}$$

which holds true provided that we take $\gamma_i=1$. This concludes the derivation of the conditional law as these $\gamma_i$'s comply with the constraint $\sum_{i=1}^{n-1} \gamma_iw_i^2=\sum_{i=1}^{n-1} w_i^2=\|w\|_2^2-w_n^2$.\medskip 


%
%

\section*{Appendix: proof of Lemma 1}

The matrix $A(m)$ is the sum of two positive definite matrices: a diagonal matrix with strictly positive entries $a_1,\ldots,a_m$ and a constant matrix with entries all set to $a_0>0$. Hence, $A(m)$ is positive definite, showing $(i)$.

Let us now compute the determinant of $A(m)$. We proceed by recursion, showing that it is true for $m+1$ whenever it holds for $m\geq 2$. It is obvious to check that it is true for $m=2$. The key point is to notice that it is enough to establish the following recursion rule :

 $$|A(m+1)|=\pi(m+1)s(m+1)=\sum_{k=0}^{m}\frac{\pi(m+1)}{a_i}+\frac{\pi(m+1)}{a_{m+1}}=a_{m+1}|A(m)|+\pi(m)\;.$$

We now apply the standard procedure for computing determinants, taking the product of each element $A(m)_{m+1,j}$ of the last row of $A(m)$ with the corresponding cofactor matrix $A(m)^{m+1,j}$ and computing the sum. Recall that the \textit{cofactor matrix associated to $A(m)_{i,j}$} is the submatrix $A(m)^{i,j}$ obtained by deleting the $i$-th row and $j$-th column of $A(m)$~\cite{Gen07}. This yields

$$|A(m+1)|=a_0\sum_{i=1}^{m}(-1)^{m+1+i}|A(m+1)^{m+1,i}|+(a_{m+1}+a_0)|A(m+1)^{m+1,m+1}|$$

where $|A(m+1)^{i,j}|$ is the minor associated to the $(i,j)$ element of $A(m)$, i.e. the determinant of the cofactor matrix $A(m+1)^{i,j}$. Interestingly, the cofactor matrices $A(m+1)^{i,j}$ take a form that is similar to $A(m)$. For instance $A(m+1)^{m+1,m+1}=A(m)$ and $A(m+1)^{m+1,m}$ is just $A(m)$ with $a_{m}\leftarrow 0$. Similarly, $A(m+1)^{m+1,1}$ is the same as $A(m)$ with $a_{1}\leftarrow 0$ provided that we shift all columns to the left, and put the last column back in first place (potentially changing the sign of the corresponding determinant), etc. More generally, for $i\in\{1,2,\ldots,m\}$, the determinant of the $(i,j)$ cofactor matrix of $A(m)$, $|A(m+1)^{i,j}|$ is exactly that of $A(m)$ with $a_i\leftarrow a_{m+1}$ if $i=j$ or that of $A(m)$ with $a_{i}\leftarrow 0$ and $a_j\leftarrow a_{m+1}$ when $j\neq i$, up to some permutations of rows and columns. In fact :

\beqn
|A(m+1)^{i,i}|&=&\sum_{k=0,k\neq i}^{m+1}\frac{\prod_{p=0}^m a_p}{a_k}\frac{a_{m+1}}{a_i}+\frac{\prod_{k=0}^m a_k}{a_i}=\frac{\pi(m+1)}{a_i}\sum_{k=0,k\neq i}^{m+1}\frac{1}{a_k}\label{eq:Minorii}\\ 
|A(m+1)^{i,j\neq i}|&=&-(-1)^{i+j}\left(\sum_{k=0,k\notin \{i,j\}}^{m+1}\frac{\pi(m+1)}{a_k}\frac{0}{a_{i}}+\frac{\pi(m)}{a_i}\frac{a_{m+1}}{a_j}\right)=-(-1)^{i+j}\frac{\pi(m+1)}{a_ia_j}\label{eq:Minorij}
\eeqn


The minor $|A(m+1)^{m+1,i}|$ when $i\neq m+1$ can be obtained from the expression of $|A(m)|$ provided that we adjust the sign and replace $a_{i}$ by 0:

$$|A(m+1)^{m+1,i}|=-(-1)^{i+m+1}\frac{\pi(m)}{a_i}~~,i\in\{1,2,\ldots,m\}$$

(recall that $A(m)$ is symmetric so that $A(m+1)^{m+1,i}=A(m+1)^{i,m+1}$). Therefore,

\beqn
|A(m+1)|&=&(a_{m+1}+a_0)|A(m)|+a_0\sum_{i=1}^{m}(-1)^{m+1+i}|A(m+1)^{m+1,i}|\nonumber\\
&=&a_{m+1}|A(m)|+a_0|A(m)|+a_0\sum_{i=1}^{m}-(-1)^{2(m+1+i)}\frac{\pi(m)}{a_i}\nonumber\\
&=&a_{m+1}|A(m)|+a_0\left(\frac{\pi(m)}{a_0} +\sum_{i=1}^{m} \frac{\pi(m)}{a_i}\right)-a_0\sum_{i=1}^{m}\frac{\pi(m)}{a_i}\nonumber\\
&=&a_{m+1}|A(m)|+\pi(m)\nonumber
\eeqn

and this recursion is equivalent to $(ii)$. \medskip

Finally, the expression of $B_{ij}(m)$ of $B(m):=(A(m))^{-1}$ are given by $1/|A(m)|$ times the adjunct matrix of $A(m)$, which is the (symmetric) cofactor matrix $C(m)$. Observe that the elements $C_{i,j}(m)$ are given by $(-1)^{i+j}M(m)_{i,j}$ where $M(m)_{i,j}$ is the minor associated to $A(m)_{i,j}$, i.e. $|A(m)^{i,j}|$. Using the minors expressions~\eqref{eq:Minorii} and~\eqref{eq:Minorij} derived above replacing $m$ by $m-1$ yields :


\beqn
B(m)_{i,i}&=&\frac{|A(m)^{i,i}|}{|A(m)|}=\frac{\sum_{k=0,k\neq j}^{m}\frac{1}{a_k}}{a_i\sum_{k=0}^m\frac{1}{a_k}}=\frac{s(m)-1/a_i}{a_is(m)}=\frac{a_is(m)-1}{a_i^2s(m)}\nonumber\\
B(m)_{i,j\neq i}&=&(-1)^{i+j}\frac{|A(m)^{i,j\neq i}|}{|A(m)|}=-\frac{\pi(m)}{a_ia_j|A(m)|}=\frac{-1}{a_ia_j\sum_{k=0}^m\frac{1}{a_k}}=\frac{-1}{a_ia_js(m)}\;.\nonumber
\eeqn




\ifdefined \MyBib
	\section{References}
  \bibliographystyle{unsrtnat}
	\bibliography{\MyBib}

\begin{thebibliography}{6}
\providecommand{\natexlab}[1]{#1}
\providecommand{\url}[1]{\texttt{#1}}
\expandafter\ifx\csname urlstyle\endcsname\relax
  \providecommand{\doi}[1]{doi: #1}\else
  \providecommand{\doi}{doi: \begingroup \urlstyle{rm}\Url}\fi

\bibitem[Li(2016)]{Li00}
D.~Li.
\newblock On default correlation: a copula function approach.
\newblock Technical report, 2016.

\bibitem[Hull and White(2004)]{Hull04}
J.~Hull and A.~White.
\newblock Valuation of a cdo and an nth-to-default cds without monte carlo
  simulation.
\newblock Technical report, 2004.

\bibitem[Andersen and Sidenius(2004)]{Ande04}
L.~Andersen and J.~Sidenius.
\newblock Extensions to the gaussian copula: random recovery and random factor
  loadings.
\newblock \emph{Journal of Credit Risk}, 1\penalty0 (1):\penalty0 29--70, 2004.

\bibitem[Vrins(2009)]{Vrins09c}
F.~Vrins.
\newblock Double t copula pricing of structured credit products - practical
  aspects of a trustworthy implementation.
\newblock \emph{Journal of Credit Risk}, 5\penalty0 (3):\penalty0 91--109,
  2009.

\bibitem[Laurent and Sestier(2016)]{Laur16}
J.-P. Laurent and M.~Sestier.
\newblock Trading book and credit risk: how fundamental is the {B}asel review ?
\newblock \emph{Journal of Banking and Finance}, 73:\penalty0 211--223, 2016.

\bibitem[Gentle(2007)]{Gen07}
J.~Gentle.
\newblock \emph{Matrix Algebra: Theory, Computations, and Applications in
  Statistics}.
\newblock Springer Texts in Statistics. Springer Texts in Statistics, 2007.

\end{thebibliography}
\fi
\end{document}